\documentclass[conference]{IEEEtran}
\IEEEoverridecommandlockouts                              
\usepackage{bbm}
\usepackage{amsfonts}
\usepackage{amssymb}
\usepackage{mathrsfs}
\usepackage{stmaryrd}
\usepackage{amsmath}
\usepackage{amssymb}
\usepackage{algorithm}
\usepackage{algorithmicx}
\usepackage{multirow}
\usepackage{tabularx}
\usepackage{latexsym}
\usepackage{times}
\usepackage{amsfonts}
\usepackage{ifthen,float}
\usepackage{color}
\usepackage{rawfonts}
\usepackage{tipa}
\usepackage{mathrsfs}
\usepackage{graphicx,color}
\usepackage{latexsym}
\usepackage{amsmath,epsfig}
\usepackage{epic,eepic}
\usepackage{cite}
\usepackage{mathrsfs}
\usepackage{float}
\usepackage{algorithm}

\floatstyle{plain}
\newfloat{twocolequfloat}{b}{zzz}
\floatname{twocolequfloat}{Equation}

\title{ Estimate the Lost Phasor Measurement Unit Data Using Alternating Direction Multipliers Method}

\author{\thanks{This work is funded by SGCC Science and Technology Program under contract No. 5455HJ160007.}
Mang Liao, \emph{Student Member, IEEE,} Di Shi, \emph{Senior Member, IEEE,} Zhe Yu, \emph{Member, IEEE,} \\Wendong Zhu, \emph{Member, IEEE,}
Zhiwei Wang, \emph{Member, IEEE,} and Yingmeng Xiang, \emph{Student Member, IEEE}
\thanks{M. Liao is with GEIRI North America, San Jose CA 95134, USA, and  the Department of Electrical and Computer Engineering,
        North Carolina State University, Raleigh, NC 27695, USA. Email: {\tt mliao@ncsu.edu}.}
\thanks{D. Shi, Z. Yu, W. Zhu, and Z. Wang are with GEIRI North America, San Jose CA 95134, USA. Email:
{\tt\{di.shi,zhe.yu,wendong.zhu,zhiwei.wang\}@geirina.net}.}
\thanks{Y. Xiang is with GEIRI North America, San Jose CA 95134, USA, and the Department of Electrical and Computer Engineering, University of Wisconsin-Milwaukee, Milwaukee, WI 53211, USA. Email: {\tt xiangy@uwm.edu}}
}

\begin{document}
\maketitle

\begin{abstract}
This paper presents a novel algorithm for recovering
missing data of phasor measurement units (PMUs). Due to the low-rank
property of PMU data, missing measurement estimation can
be formulated as a low-rank matrix-completion problem. Based
on maximum-margin matrix factorization, we propose an efficient
algorithm based on alternating direction method of multipliers (ADMM)
for solving the matrix completion problem. Comparing to existing
approaches, the proposed ADMM based algorithm does not need
to estimate the rank of the target data matrix and provides
better performance in computation complexity.
In addition, we consider the case of measurements missing from all PMU
channels and provide a strategy
of reshaping the matrix which contains the received PMU data for
estimation. Numerical results using PMU measurements
from IEEE 68-bus power system model illustrate the
effectiveness and efficiency of the proposed approaches.
\end{abstract}


%


\maketitle

\pagenumbering{arabic}

\section{Introduction}
The wide-area measurement system (WAMS) using phasor
measurement units (PMUs) has been regarded as one of the
key enabling technologies in monitoring, control, and protections
of the next-generation power grids \cite{Phadke2008}. With continuous
increase in PMU deployment and the resulting explosion in
data volume, the design and deployment of an efficient wide area
communication and computing infrastructure, especially
from the point of view of resilience against a large number of
missing data, is evolving as one of the greatest challenges
to the power system and IT communities.
With thousands of networked PMUs
being scheduled to be installed in the United States by 2020,
exchange of synchrophasor data between balancing authorities
for any type of wide-area control will involve an enormous number of data flow in real-time per event, thereby
opening up a wide spectrum of probabilities of data
losses and data quality degradations in an unpredictable way. Data missing makes the system unobservable, degrades the
performance of the state estimates, and weakens the security
and stability of the system. Therefore, recovering missing PMU measurements has become a significant and inevitable
problem in power systems.

PMU data can be structured as a matrix with each column and row representing the measurements of one channel and sample instant, respectively. Since large amounts of PMU data exhibit heavily correlated property \cite{Chen2013, Dahal2012, Gao2014}, the matrix is approximately low-rank, and the problem of recovering the missing PMU data can be formulated as a {\it low-rank matrix-completion} problem.
Studies on matrix completion algorithms are extensive, including atomic decomposition of minimum rank approximation (ADMiRa) \cite{Lee2010}, singular value projection (SVP) \cite{Jain2010}, information cascading matrix completion (ICMC) \cite{Meka2009}, among which nuclear-norm-regularized matrix approximation \cite{Cand2s2010,Mazumder2010,Cai2010,Ma2011} and maximum-margin matrix factorization (MMMF) \cite{Rennie2005} are widely adapted. Using nuclear-norm-regularized matrix approximation, a singular value threshold has to be designed which influences the estimate accuracy. Developing the nuclear-norm-regularized matrix approximation, an alternating direction method (ADM) is provided for solving the matrix completion problem \cite{Chen2012,Xu2016}.
However, the calculation of the singular value decomposition (SVD) in ADM approach increases the computational time and complexity. Based on MMMF, Jain et al. \cite{Jain2013} and Hardt \cite{Hardt2014} proposed alternating least squares (ALS) schemes for solving the matrix-completion problem. Further, {\it softImpute-ALS} is provided for reducing the computational complexity \cite{Hastie2015}. Gao et al. applied the MMMF approach on recovering the missing PMU data \cite{Gao2016} firstly.
Most of the existing approaches rely on an estimation of
the rank $r$ of the data matrix, which is typically unavailable and time variant in practice. Inaccurate estimation of $r$ introduces modelling errors in the matrix completion problem. The computational complexity is lower with a smaller $r$. On the other
hand $r$ cannot be too small for estimation accuracy.
Therefore, design of an adaptive and scalable online algorithm
of PMU data recovery is an open challenge.

Motivated by these insights, we develop an algorithm that can recover the missing PMU measurement with low computational complexity and less operating time.
The fundamental set-up for this optimization was based on  MMMF and alternating direction method of multipliers (ADMM) \cite{Glwinski&Marroco1975,Gabay&Mercier1976,Glowinski1989}.
 Firstly, the observed PMU data is structured as a matrix ${{\boldsymbol M} \in \mathbb{R}^{n_1 \times n_2}}$ whose columns and rows represent the measurements from one channel and the same sampling instant, respectively. Then we formulate the data recovery as an optimization problem in which we minimize the rank of the estimated matrix $\hat{\boldsymbol X}$ while keeping elements in $\hat{\boldsymbol X}$ the same as the corresponding ones in $\boldsymbol M$ if they are present.
 An ADMM algorithm is proposed to solve the optimization problem in an iterative way. In the update equations there is no matrix inverse computation, which immensely reduces the computational complexity. In addition, it is not necessary to estimate the rank of the original data matrix ${\boldsymbol X}$ without missing elements, which significantly cuts down the influence of the uncertain factor into the performance.
Furthermore, we consider the case of missing data from all
PMU channels. In this case, all elements in one row of the observed matrix $\boldsymbol M$ are missing. One efficient algorithm is presented to
reshape the observed matrix, and the lost data from all the
channels can be recovered
using ADMM approach.
We illustrate the results using simulations of the IEEE 68-bus system model.

\section{Problem Formulation}
Persistent model is one simple and traditional method to recover the missing PMU data. It utilizes the temporal
correlation of the PMU measurements to recover the lost
data in one channel. However, if in the disturbance scenario the measurements in the same channel are missing during a long time, the recovery with persistent model is not an advisable choise.
In this section, we process a spatial-temporal blocks of PMU data, present the low-rank property of PMU
data, and formulate the data recovery as a matrix completion
problem.

\subsection{Low-rank property of PMU measurements}
Denote ${{\boldsymbol X} \in \mathbb{R}^{n_1 \times n_2}}$ as the PMU measurement matrix
without data missing. Each column and row correspond to a sequence of measurements of one PMU channel, and the PMU measurments at the same sampling instant, respectively. Due
to the noise, all the singular values of $\boldsymbol X$ are larger than zero. An approximating rank approach, referred to Frobenius norm proportion \cite{Zhang2004}, is stated as follows.
\begin{align}\label{apro_rank}
\frac{\sqrt{\sigma_1^2 + \sigma_2^2 + \ldots \sigma_r^2}}{\sqrt{\sigma_1^2 + \sigma_2^2 + \ldots \sigma_r^2 + \ldots + \sigma_l^2}} \ge \beta,
\end{align}
where ${\sigma_1> \sigma_2> \ldots > \sigma_l}$ are the singular values of the matrix and $\beta$, $0 < \beta \le 1$, is the proportion factor. $r$ in \eqref{apro_rank} denotes the approximate rank of the matrix.
Since the PMU measurements of voltage or current phasors or magnitudes from different lines or buses are strongly correlated, the approximate rank of $\boldsymbol X$ is much smaller than ${\rm min}\{n_1, n_2\}$ \cite{Chen2013, Dahal2012, Gao2014}.
Due to the low-rank property of PMU data, missing
PMU measurement estimation can be converted into a low-rank
matrix completion problem.

%
%

\subsection{ An ADMM based approach for PMU data estimation}\label{ADMM_sec}

Let ${{\boldsymbol M} \in \mathbb{R}^{n_1 \times n_2}}$ and ${\hat{\boldsymbol X}} \in \mathbb{R}^{n_1 \times n_2}$ denote the observed PMU measurements with missing data and the recovered matrix, respectively. Since $\hat{\boldsymbol X}$ should be a low-rank matrix, the matrix completion problem is formulated as
follows:
\begin{equation}\label{ori_object1}
\begin{array}{ll}
\min\limits_{{\hat{\boldsymbol X}} \in \mathbb{R}^{n_1 \times n_2}} & {\rm rank}({\hat{\boldsymbol X}})\\
\mbox{subject to} & ({\hat{\boldsymbol X}} - {\boldsymbol M})\odot {\boldsymbol I}_s = {\boldsymbol 0},
\end{array}
\end{equation}
where $\odot$ denotes the Hadamard product, \emph{i.e.}, ${[{\boldsymbol Y}_1 \odot {\boldsymbol Y}_2]_{ij} = [{\boldsymbol Y}_1]_{ij}[{\boldsymbol Y}_2]_{ij}}$.
${\boldsymbol I}_s$ is the {\it structural identity} with its $ij^{th}$ entry defined as
\begin{align}
{\left[ {{{\boldsymbol I}_s}} \right]_{ij}} = \left\{ {\begin{array}{*{20}{l}}
{1,}&{if~{[{\boldsymbol M}]_{ij}}~is~observed~data;}\\
{0,}&{if~{[{\boldsymbol M}]_{ij}} ~is~missing~data.}
\end{array}} \right.
\end{align}
Unfortunately, \eqref{ori_object1} is NP hard to solve, and can be relaxed to a tractable optimization problem \cite{PhD}:
\begin{equation}\label{ori_object}
\begin{array}{ll}
 \min \limits_{{\hat{\boldsymbol X}} \in \mathbb{R}^{n_1 \times n_2}}& ||{\hat{\boldsymbol X}}||_{*}\\
\mbox{subject to} & ({\hat{\boldsymbol X}} - {\boldsymbol M})\odot {\boldsymbol I}_s = {\boldsymbol 0},
\end{array}
\end{equation}
where the nuclear norm $||{\hat{\boldsymbol X}}||_*$ is the sum of the singular values of ${\hat{\boldsymbol X}}$.

Using MMMF to further change the optimization problem \eqref{ori_object}, let ${{\hat{\boldsymbol X}} = {\boldsymbol A}^{\rm T} {\boldsymbol B}}$, in which ${{\boldsymbol A} \in \mathbb{R}^{n_2 \times n_1}}$ and ${{\boldsymbol B} \in \mathbb{R}^{n_2 \times n_2}}$. Without loss of generality, we assume ${n_1 > n_2}$. Since
$||{\hat{\boldsymbol X}}||_*$ is equivalent to $\mathop{\min}\limits_{{\boldsymbol A}, {\boldsymbol B}}\frac{1}{2}(||\boldsymbol A||_F^2 + ||{\boldsymbol B}||_F^2) $
with Frobenius norm $||.||_F$ \cite{Rennie2005}, the optimization function is equivalent to
\begin{equation}\label{objective}
\begin{array}{ll}
\min \limits_{{\boldsymbol A}, {\boldsymbol B} }& \frac{1}{2}(||{\boldsymbol A}||_{F}^2 + ||\boldsymbol B||_F^2)\\
\mbox{subject to} & ({\boldsymbol A}^{\rm T}{\boldsymbol B} - {\boldsymbol M})\odot {\boldsymbol I}_s = {\boldsymbol 0}.
\end{array}
\end{equation}
In the previous work \cite{Jain2013, Hardt2014}, people estimated the rank $r$ of ${\hat{\boldsymbol X}}$, set ${{\boldsymbol A} \in \mathbb{R}^{r \times n_1}}$ and ${{\boldsymbol B} \in \mathbb{R}^{r \times n_2}}$, and applied ALS to solve \eqref{objective}. The computational complexity is ${\mathcal{O}((n_1+ n_2)r^{3})}$. If ${r = {\rm min}\{n_1, n_2\}}$, the computational complexity is $\mathcal{O}((n_1 + n_2)({\rm min}\{n_1, n_2\})^3)$, which is a biquadrate function of ${\rm min}\{n_1, n_2\}$. With smaller $r$ the computational complexity is reduced.
However, the value of $r$ cannot be too small to guarantee the estimation accuracy. For reducing the influence of the uncertain factor into the performance, we set the sizes of matrices ${\boldsymbol A\in\mathbb{R}^{n_2 \times n_1}}$ and ${\boldsymbol B\in\mathbb{R}^{n_2 \times n_2}}$ only depend on the size of observed matrix $\boldsymbol M$. In addition, we apply the ADMM method to solve \eqref{objective} in an iterative way using the Lagrangian multiplier approach.

The augmented Lagrangian for \eqref{objective} can be formulated as
\begin{align}\label{org_L}
\mathcal{L} =& \frac{1}{2}(||{\boldsymbol A}||_{F}^2 + ||\boldsymbol B||_F^2) + {\rm trace}({\boldsymbol w}^{\rm T}(({\boldsymbol A}^{\rm T}{\boldsymbol B} - {\boldsymbol M})\odot {\boldsymbol I}_s ))+ \nonumber\\
& \frac{\rho}{2}||({\boldsymbol A}^{\rm T}{\boldsymbol B} - {\boldsymbol M})\odot {\boldsymbol I}_s ||_F^2,
\end{align}
where $\boldsymbol A$ and ${\boldsymbol B}$ are the matrices of the primal variables, ${\boldsymbol w}$ is the matrix of the dual variables or the Lagrange multipliers associated with \eqref{objective}, and $\rho >0$ denotes a penalty weight.

After some algebraic, the augmented Lagrangian can be rewritten as
\begin{equation}\label{changed_L}
\begin{array}{l}
\mathcal{L} = \frac{1}{2}(||{\boldsymbol A}||_{F}^2 + ||\boldsymbol B||_F^2) + {\rm trace}(({\boldsymbol w}\odot {\boldsymbol I}_s)^{\rm T}({\boldsymbol A}^{\rm T}{\boldsymbol B} - {\boldsymbol M}))+\\
\mathrel{\phantom{\mathcal{L} =}}  \frac{\rho}{2}{\rm trace}((({\boldsymbol A}^{\rm T}{\boldsymbol B} - {\boldsymbol M})\odot {\boldsymbol I}_s )^{\rm T}({\boldsymbol A}^{\rm T}{\boldsymbol B} - {\boldsymbol M})).
\end{array}
\end{equation}
The gradients of the augmented Largrangian $\mathcal{L}$ in \eqref{changed_L} with respect to $\boldsymbol A$ and ${\boldsymbol B}$ are respectively given by
\begin{align}
\frac{{\partial {\mathcal L}}}{{\partial {\boldsymbol A}}} &= {\boldsymbol A} + {\boldsymbol B}({\boldsymbol w}\odot {\boldsymbol I}_s)^{\rm T}+\rho{\boldsymbol B}(({\boldsymbol A}^{\rm T}{\boldsymbol B} - {\boldsymbol M})\odot {\boldsymbol I}_s )^{\rm T},\nonumber\\
\frac{{\partial {\mathcal L}}}{{\partial {\boldsymbol B}}} &= {\boldsymbol B} +  {\boldsymbol A}({\boldsymbol w}\odot {\boldsymbol I}_s) + \rho{\boldsymbol A}({\boldsymbol A}^{\rm T}{\boldsymbol B} - {\boldsymbol M})\odot {\boldsymbol I}_s.
\end{align}
Given the derivation, the ADMM algorithm for solving the optimal problem \eqref{objective} is
illustrated in Algorithm \ref{ADMM}.
\begin{algorithm}
\caption{ADMM algorithm for PMU data estimation}
\begin{algorithmic}
\State Initialize ${\boldsymbol A}^0$, ${\boldsymbol B}^0$, ${\boldsymbol w}^0$, and $k=0$, and determine the value of $\rho$, $\epsilon$ and $k_{\max}$.
\State{\bf Do:}
\begin{small}
\begin{equation}
\hspace{-1em}
\begin{array}{l}
{\boldsymbol A}^{k+1} = -{\boldsymbol B}^k({\boldsymbol w}^k\odot {\boldsymbol I}_s)^{\rm T}-\rho{\boldsymbol B}^k((({\boldsymbol A}^k)^{\rm T}{\boldsymbol B}^k - {\boldsymbol M})\odot {\boldsymbol I}_s )^{\rm T}, \\
{\boldsymbol B}^{k+1} = -  {\boldsymbol A}^{k+1}({\boldsymbol w}^k\odot {\boldsymbol I}_s) - \rho{\boldsymbol A}^{k+1}((({\boldsymbol A}^{k+1})^{\rm T}{\boldsymbol B}^k - {\boldsymbol M})\odot {\boldsymbol I}_s), \\
{\boldsymbol w}^{k+1} = {\boldsymbol w}^k + \rho(({\boldsymbol A}^{k+1})^{\rm T}{\boldsymbol B}^{k+1} - {\boldsymbol M})\odot {\boldsymbol I}_s,\\
k=k+1.
\end{array}
\end{equation}
\end{small}
\State{\bf until:}\\
 The stopping criterion ${||({\boldsymbol A}^{k+1})^{\rm T}{\boldsymbol B}^{k+1} - ({\boldsymbol A}^{k})^{\rm T}{\boldsymbol B}^{k}||<\epsilon}$ is reached or ${k>k_{\max}}$.
\end{algorithmic}\label{ADMM}
\end{algorithm}


The updates in Algorithm \ref{ADMM} requires no matrix inverse, and the computational complexity is $\mathcal{O}(n_1 n_2 {\rm min}\{n_1, n_2\})$, which is a quadratic function of ${\rm min}\{n_1, n_2\}$. In addition, it is not necessary to estimate the rank of matrix ${\hat{\boldsymbol X}}$, which reduces the influence of uncertain factor into the performance. Penalty
weight $\rho$ denotes the step size of the dual variable
update. In general, large $\rho$ results in fast convergent rate.

Compared to approaches like interpolations and persistent models, ADMM algorithm utilizes the spatial and temporal correlations of PMU data to improve accuracy. In the persistent model, it replaces the missing data by the previous available data point. The persistent method recovers the lost data only based on temporal correlation. If the data from one channel are missing during a long time, and if there exists a dynamic in the time, then the estimation using persistent method doubtlessly is a nightmare. On the other hand, based on the spatial correlation, the missing data can be recovered using ADMM. We will compare the estimates using ADMM and persistent model with IEEE 68-bus power system
simulation in Subsection \ref{Comparison}.

\subsection{Special case: missing data from all the channels}\label{reshape}

The power system often suffers natural and artificial disturbances during operation. It is possible that the data from all the channels are missing simultaneously under communication failure. In this case, no existing algorithms can recover the missing data.
For solving this problem, the observed matrix $\boldsymbol M$ has to be reshaped to avoid some of its rows missing.
Our goal is that the proportion of missing elements in one row of the reshaped observed matrix $\boldsymbol M$ is as small as possible. Meanwhile the corresponding reshaped recovery matrix $\hat{\boldsymbol X}$ is still low-rank. 

We provide an alternative method, called {\it cut-column reshaping method} (CCRM), for reshaping the observed matrix.
Using CCRM each column with $n_1$ length is separated into $n^*$ shorter columns with a length of $\frac{n_1}{n^*}$. Thus, the $n_1$-by-$n_2$ matrix is reshaped to a $\frac{n_1}{n^*}$-by-$n_2n^*$ matrix, and the original column correlation is held. The length of the new column should be larger than the row length of the original matrix, \emph{i.e.}, $\frac{n_1}{n^*} > n_2$. $n^*$ also satisfies that $\left\lceil\frac{n_1}{n^*+1}\right\rceil<n_2$, where $\left\lceil x \right\rceil$ denotes the smallest integer number which is larger than $x$.
Thus the numbers of rows and columns of reshaped matrix are both larger than $n_2$. Due to the size, the rank of $\boldsymbol M$ is no more than ${\rm min}\{n_1, n_2\}$. Using CCRM the rank of reshaped matrix ${\tilde{\boldsymbol M}}$ will not be reduced by the new size. In addition, with holding the column correlation, CCRM minimizes the proportion of zero elements in one row of reshaped matrix.

Consider a simple example to illustrate the reshaping method.
A $6 $-by-$2$ matrix $\boldsymbol M$ can be expressed as:
\begin{equation}
\label{ex_M}
\begin{array}{lll}
{\boldsymbol M} &= \left[{\begin{array}{*{20}{c}}
{\boldsymbol m}_{1} &{\boldsymbol m}_{2}
\end{array}} \right] \\
&= \left[ {\begin{array}{*{20}{c}}
{m}_{11} & {m}_{21} &  m_{31}  & m_{41} &  \star &  m_{61}\\
{m}_{12} & {m}_{22} &  m_{32}  & m_{42} &  \star &  m_{62}
\end{array}} \right]^{\rm T}
\end{array}
\end{equation}
whose fifth row is missing. Using CCRM with ${n^* = 3}$ and matrix ${\boldsymbol M}$ is reshaped into a $2$-by-$6$ matrix:
\begin{align}
\begin{array}{ll}
\tilde{\boldsymbol M} &= \left[ {\begin{array}{*{20}{c}}
{\tilde{\boldsymbol m}}_{1}&{\tilde{\boldsymbol m}}_{2}
&{\tilde{\boldsymbol m}}_{3}&{\tilde{\boldsymbol m}}_{4}
&{\tilde{\boldsymbol m}}_{5}&{\tilde{\boldsymbol m}}_{6}
\end{array}} \right] \\
&= \left[ {\begin{array}{*{20}{c}}
m_{11} & m_{31} & \star & m_{12} & m_{32} & \star \\
m_{21} & m_{41} & m_{61} & m_{22} & m_{42} & m_{62}
\end{array}} \right].
\end{array}
\end{align}
Now for each column and row, not all measurements are missing. If ${\boldsymbol m}_1$ and ${\boldsymbol m}_2$ are strongly correlated, ${\tilde{\boldsymbol m}}_{1}$ and ${\tilde{\boldsymbol m}}_{4}$, ${\tilde{\boldsymbol m}}_{2}$ and ${\tilde{\boldsymbol m}}_{5}$, and ${\tilde{\boldsymbol m}}_{3}$ and ${\tilde{\boldsymbol m}}_{6}$ are strongly correlated in pairs. The ranks of matrice ${\boldsymbol M}$ and ${\tilde{\boldsymbol M}}$ are both no more than $2$. The proportion of missing elements to the first row is $\frac{1}{3}$; while it is $1$ to the fifth row of ${\boldsymbol M}$.
CCRM is illustrated in Algorithm \ref{CCRM}.
\begin{algorithm}
\caption{Cut-Column Reshaping Method}
\begin{algorithmic}
\State(1) Check whether any row of the observed $n_1$-by-$n_2$ matrix ${\boldsymbol M}$ owns all missing elements.
\State(2) If yes, let $n^*$ be the maximum divisor of $n_1$, which satisfies $\frac{n_1}{n^*} >n_2$.
\State(3) Separate each column of $\boldsymbol M$ into $n^*$ shorter columns with $\frac{n_1}{n^*}$ length. The original $n_1$-by-$n_2$ matrix is reshaped into a $\frac{n_1}{n^*}$-by-$n_2 n^*$ matrix.
\end{algorithmic}\label{CCRM}
\end{algorithm}
The missing PMU measurements from all the channels can be recovered using ADMM in Algorithm \ref{ADMM} after reshaped matrix $\boldsymbol M$ using CCRM in Algorithm \ref{CCRM}.
Notice that if all the elements in one column of the reshaped observed matrix are missing, they cannot be recovered using ADMM. The recovery accuracy using ADMM will be declined sharply, if the measurements in one channel are missing more than $\frac{n_1}{2}$ successive sampling instants. With less lost data, the recovery accuracy will be enhanced.

\section{Simulation results}
The IEEE 68-bus system is used to carry out the simulation to verify the proposals. We build up a PMU measurement matrix whose column and row corresponding to a sequence voltage phasors on $86$ lines and the sampling instants, respectively. The simulated measurements are obtained using the power systems
toolbox (PST) nonlinear dynamics simulation routine $\it s\_simu$
and the data file {\it data16m.m} \cite{Chow1992}. A three-phase fault is imposed at the line connecting buses $1$ and $2$. The fault starts at $t = 0.1$s, and clears on bus $1$ at $t = 0.15$s and on bus $2$ at $t = 0.20$s.
For approaching to the true measurements, we add white Gaussian noise ($\mathcal{N}(0, 0.001)$) into the PMU data.
The measurements are observed during $60$s and
there are $30$ samples in one second. The $1800$-by-$86$ matrix ${\boldsymbol X}$ is with no
missing measurements and its approximate rank is $1$ with
$\beta = 0.995$ in \eqref{apro_rank}.
To test the recovery accuracy of the presented ADMM algorithm,
some observed data in $\boldsymbol X$ is set to be lost. Since
the PMU data are missing arbitrary and unpredictable, in this
paper we consider two cases of missing data: (1) Missing data
randomly. The delivery of PMU measurements from multiple
remote locations of power grids to monitoring centers can
result in the random unavailability of PMU measurements;
(2) Missing data in all channels simultaneously. The
transform link malfunctions may  result in data missing in all channels.
We choose the penalty weight ${\rho = 0.00075}$ using ADMM, and the dual parameter ${\lambda = 1.5}$ and the estimated rank of filled completion matrix ${r = 20}$ using ALS for comparison.
In the paper, the computational
time is obtained by operating Matlab programming.

\subsection{Case 1: Missing data randomly}

In this case, we assume an independent and identical
distribution (i.i.d) of the missing rate. For each data point, with
a probability the measurement is missing
and set to zero in $\boldsymbol M$ artificially. Notice that it is different from the data which is equal to zero. If the actual data is zero, the corresponding element in ${\boldsymbol I}_s$ is equal to $1$. While if the data is missing, the corresponding element in ${\boldsymbol I}_s$ is equal to $0$.
%
Table \ref{suma} compares some properties of ALS and ADMM in Case 1.
\begin{table}[h]
\centering
\caption{Comparison of ALS and ADMM for the recovery}
\begin{tabular}{|c|c|c|c|}
\hline
& \# iterations & time & Sensitivity of parameters\\
\hline
ALS & $\approx 50$ & $>7s$   & Less stringent \\
\hline
ADMM & $\approx 100$ &  $<1s$  & More stringent \\
\hline
\end{tabular}
\label{suma}
\end{table}
Though the number of convergence iterations using ADMM is larger than the one using ALS, the computational time using ADMM is less than $1s$, which is much smaller than using ALS. 

\begin{figure}[!t]
\centering
\includegraphics[width=8.3cm,height=5.3cm]{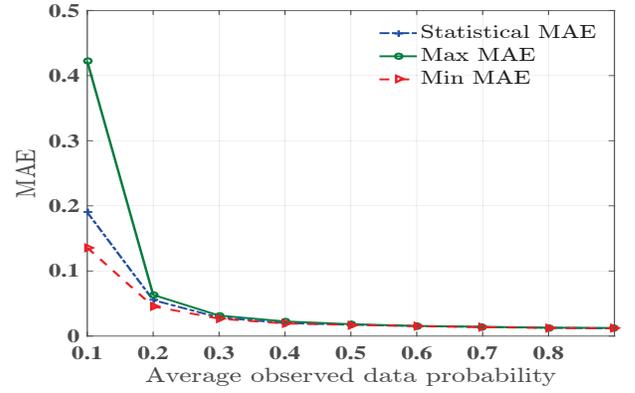}
\caption{Case 1: MAEs using ADMM and ALS against different observed data probabilities, respectively.}
\label{MAE_1}
\end{figure}

Fig. \ref{MAE_1} shows the statistic, maximum, and minimum values of mean absolute errors MAEs $\frac{\sum\limits_{ij:[{\boldsymbol I}_s]_{ij}=0}{|[{\hat{\boldsymbol X}}]_{ij} - [{\boldsymbol X}]_{ij}|}}{\sum\limits_{ij:[{\boldsymbol I}_s]_{ij}=0}[{\boldsymbol I}_s]_{ij}}$ using ADMM and ALS with different observed data probabilities, respectively. The observed data probability denotes the likelihood of the observed data occurrence, i.e., $\frac{\sum[{\boldsymbol I}]_{ij}}{n_1n_2}$.The statistic, maximum, and minimum values of MAEs are obtained by Monte Carlo method with $500$ independent times. With larger probability of observed measurements, MAE becomes smaller. The
statistic values of MAEs using ADMM and ALS are close with each observed data probability. The difference between the maximum and minimum values of MAEs using ADMM is larger than the one using ALS. 

\subsection{Case 2: Missing data in all channels}

In this case, one row of data in matrix $\boldsymbol M$ is lost. The $1800$-by-$86$ matrix $\boldsymbol M$ which contains  voltage phasor measurements can be treated as $1800$ sub-matrices with a size of  $1$-by-$86$. The observed data
probability denotes the proportion of the observed sub-matrices
to the total ones.
For recovering the missing data in one row, firstly we reshape the observed matrix using CCRM. Since the rank of the orginal matrix is no more than $86$, the number of the rows and columns of the reshaped matrix should be more than $86$ for avoiding reducing the rank artificially. In addition, with holding the column correlation, one purpose of reshaping is
minimizing the proportion of zero elements in one row of the
reshaped matrix. Thus using CCRM, the
original $1800$-by-$86$ matrix $\boldsymbol M$ is reshaped to a $90$-by-$1720$ matrix ${\tilde{\boldsymbol M}}$ with
${n^* = 20}$. With ${\beta = 0.995}$ in \eqref{apro_rank}, the approximate rank of the reshaped observed matrix $\tilde{{\boldsymbol X}}$ is $1$. Since the size of the transposed reshaped observed matrix is similar to the observed matrix, the computational time using ADMM and ALS is similar to the results in Table \ref{suma}, respectively.



\begin{figure}[!t]
\centering
\includegraphics[width=8.5cm,height=5.5cm]{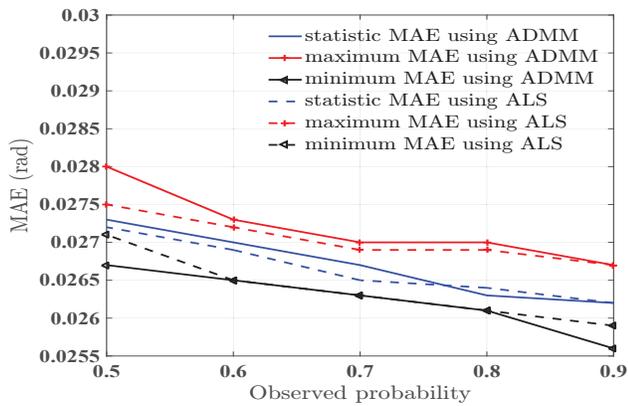}
\caption{Case 2: MAEs using ADMM and ALS against different observed data probabilities, respectively.}
\label{MAE_2}
\end{figure}

Fig. \ref{MAE_2} shows the statistic, maximum, and minimum values of MAEs using ADMM and ALS with different observed data probabilities, respectively. The
statistic values of MAEs using ADMM and ALS are still close. Compared with Case 1, the MAEs using both ADMM and ALS are larger. Though approximate rank of reshaped matrix $\tilde {\boldsymbol X}$ is still 1,
the minimum singular value becomes larger, whose influence
into the recovery accuracy cannot be ignored. If the observed matrix $\boldsymbol X$ is not reshaped, the missing row cannot be recovered using neither ADMM nor ALS, and the MAEs with different observed data probabilities are all around $0.278$.

\subsection{Comparison among ADMM, ALS, and persistent model approaches}\label{Comparison}

\begin{figure}[!t]
\centering
\includegraphics[width=8.5cm,height=5.5cm]{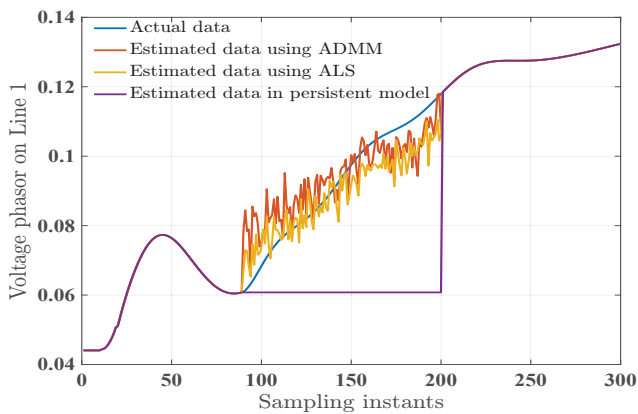}
\caption{Comparison of the estimated measurements using ADMM, ALS, and persistent model. The blue line shows the actual measurements.}
\label{comp_A_A_p2}
\end{figure}

In the persistent model, it replaces the missing data at the $t^{th}$ sampling instant with the data at the $(t-1)^{th}$ if it is available.
Only based on the temporal correlation of the PMU data, in a disturbed scenario the data which are lost during several successive sampling instants cannot be recovered successfully using persistent method. In this subsection, we let the data be lost from the $90^{th}$ sampling instant to the $200^{th}$ sampling instant on $9$ lines. Fig. \ref{comp_A_A_p2} shows the estimated measurements using ADMM, ALS, and persistent methods from sampling instant $1$ to $300$ on Line 1. The blue line shows the values of actual measurements. Using both ADMM and ALS approaches, estimated measurements are close to the actual one. While using the persistent model, the estimate deviates from the actual data due to the dynamics in the measurements.

\section{Conclusion}

In this paper, we presented ADMM
algorithm for missing PMU measurement recovery.
We illustrated our
results  with noisy measurements from the IEEE 68-bus power system model. Compared with the ALS algorithm, the computational complexity and operating time are much smaller using the ADMM algorithm. In addition, the ADMM algorithm avoids to estimate the rank of filled completion matrix, which reduces the influence of the uncertain factor into the performance. We also consider the case of missing data in all the channels simultaneously and provide one approach to reshape the observed matrix for the recovery.
Our future work in this area will include recovering continuous several rows of the observed matrix with all missing elements and testing the proposal using actual PMU data.


\begin{thebibliography}{99}


\bibitem{Phadke2008}
A. G. Phadke and J. S. Thorp, ``Synchronized Phasor Measurements and Their Applications," New York, Springer, 2008.

\bibitem{Chen2013}
Y. Chen, L. Xie, and P. Kumar, ``Dimensionality Reduction and Early Event Detection Using Online Synchrophasor Data," in {\it Proc. IEEE Power and Energy Society General Meeting}, pp. 1-5, 2013.

\bibitem{Dahal2012}
N. Dahal, R. L. King, and V. Madani, ``Online Dimension Reducion of Synchrophasor Data," in {\it Proc. IEEE PES Transmission and Distribution Conf. Expo. (T\& D)}, pp. 1-7, 2012.

\bibitem{Gao2014}
P. Gao, M. Wang, S. Ghiocel, and J. H. Chow, ``Modeless Reconstruction of Missing Synchrophasor Measurements," in {\it Proc. IEEE PES General Meeting}, pp. 1-5, 2014.


\bibitem{Lee2010}
 K. Lee and Y. Bresler, ``ADMiRA: Atomic Decomposition for Minimum Rank Approximation," {\it IEEE Trans. Inf. Theory}, vol. 56, no. 9, pp. 4402-4416, 2010.
 

 \bibitem{Jain2010}
 P. Jain, R. Meka, and I. S. Dhillon, ``Guaranteed Rank Minimization via Singular Value Projection," {\it in Adv. Neural Inf. Process. Syst.}, pp. 937-945, 2010.

 \bibitem{Meka2009}
R. Meka, P. Jain, and I. S. Dhillon, ``Matrix Completion from
Power-Law Distributed Samples," {\it in Adv. Neural Inf. Process. Syst.}, pp. 1258-1266, 2009.

\bibitem{Cand2s2010}
E. J. Candes and T. Tao, ``The Power of Convex Relaxation: Near-Optimal Matrix Completion," {\it IEEE Trans. on Information Theory}, vol. 56, no. 5, pp. 2053-2080, 2010.

\bibitem{Mazumder2010}
R. Mazumder, H. Trevor, and T. Robert, ``Spectral Regularization Algorithms for Learning Large Incomplete Matrices," {\it Journal of machine learning research}, vol. 11, pp. 2287-2322, 2010.

\bibitem{Cai2010}
J. Cai, E. J. Candes, and Z. Shen, ``A Singular Value Thresholding Algorithm for Matrix Completion," {\it SIAM Journal on Optimization}, vol. 20, no. 4, pp. 1956-1982, 2010.

\bibitem{Ma2011}
S. Ma, D. Goldfarb, and L. Chen, ``Fixed Point and Bregman Iterative Methods for Matrix Rank Minimization," {\it Mathematical Programming}, vol. 128., no. 1, pp. 321-353, 2011.


\bibitem{Rennie2005}
J. D. Rennie and S. Nathan, ``Fast Maximum Margin Matrix Factorization for Collaborative Prediction," in {\it Proc. the 22nd international conference on Machine learning(ACM)}, pp. 713-719, 2005.


\bibitem{Chen2012}
C. Chen, B. He, and X. Yuan, ``Matrix Completion via an Alternating Direction Method," {\it IMA Journal of Numerical Analysis}, vol. 32, no. 1, pp. 227-245, 2012.


\bibitem{Xu2016}
F. Xu, and P. Pan, ``A New Algorithm for Positive Semidefinite Matrix Completion," {\it Journal of Applied Mathematics}, vol. 3, pp. 1-5, 2016.


\bibitem{Jain2013}
P. Jain, N. Praneeth, and S. Sujay, ``Low-Rank Matrix Completion Using Alternating Minimization." in {\it Proc. the forty-fifth annual ACM symposium on Theory of computing (ACM)}, 2013.


\bibitem{Hardt2014}
M. Hardt, ``Understanding Alternating Minimization for Matrix Completion," {\it Foundations of Computer Science (FOCS), 2014 IEEE 55th Annual Symposium}, 2014.


\bibitem{Hastie2015}
T. Hastie, R. Mazumder, J. D. Lee, and R. Zadeh, ``Matrix Completion and Low-Rank SVD via Fast Alternating Least Squares," {\it J. Mach. Learn. Res}, vol. 16, no. 1, pp. 3367-3402, 2015.

\bibitem{Gao2016}
P. Gao, M. Wang, S. G. Ghiocel, J. H. Chow, B. Fardanesh, and G. Stefopoulos, ``Missing Data Recovery by Exploiting Low-Dimensionality in Power System Synchrophasor Measurements," {\it IEEE Trans. on Power Systems}, vol. 31, no. 2, pp. 1006-1013, 2016.



\bibitem{Glwinski&Marroco1975}
R. Glowinski and A. Marroco, ``Sur l'approximation, par \`el\`ements finis d'ordre un, et la r\`esolution, par p\`enalisation-dualit\`e d'une classe de probl\`emes de Dirichlet non lin\`eaires.'' {\it ESAIM: Mathematical Modelling and Numerical Analysis - Mod¨¦lisation Math¨¦matique et Analyse} Num¨¦rique 9.R2 pp. 41-76, 1975.



\bibitem{Gabay&Mercier1976}
D. Gabay and B. Mercier, ``A Dual Algorithm for the Solution of Nonlinear Variational Problems via Finite Element Approximation''. {\it Computers \& Mathematics with Applications}, vol. 2, no. 1, pp. 17-40, 1976. 


\bibitem{Glowinski1989}
R. Glowinski and L. T. Patrick, ``Augmented Lagrangian and operator-splitting methods in nonlinear mechanics," Society for Industrial and Applied Mathematics, 1989.


\bibitem{Zhang2004}

X. Zhang, ``Matrix Analysis and Applications," {Tsinghua and Springer Publishing house}, Beijing, pp. 71-100, 2004.


\bibitem{PhD}
M. Fazel, ``Matrix Rank Minimization with Applications," Ph. D. dissertation, Stanford Univ., Stanford, CA, 2002. 

\bibitem{Chow1992}
J. H. Chow and K. W. Cheung, ``A Toolbox for Power System
Dynamics and Control Engineering Education and Research," {\it IEEE
Transactions on Power Systems}, vol. 7, no. 4, pp. 1559-1564, 1992.
\end{thebibliography}
\end{document}